\documentclass{amsart}

\DeclareMathOperator{\Symm}{Symm}
\DeclareMathOperator{\diag}{diag}
\DeclareMathOperator{\Tr}{Tr}

\newtheorem{theorem}{Theorem}[section]
\newtheorem{lemma}[theorem]{Lemma}

\theoremstyle{definition}
\newtheorem{definition}[theorem]{Definition}

\theoremstyle{remark}
\newtheorem{remark}[theorem]{Remark}

\numberwithin{equation}{section}
\usepackage{xcolor}

\begin{document}

\title[Geodesics in the Deep Linear Network]{Geodesics in the Deep Linear Network}


\author{Alan Chen}
\address{Division of Applied Mathematics, Brown University, Providence, RI 02912}
\email{alan\_chen1@alumni.brown.edu}




\begin{abstract}
We derive a general system of ordinary differential equations and associated explicit solutions in a special case for geodesics between full rank matrices in the deep linear network geometry. In the process, we find horizontal straight lines in the invariant balanced manifold that remain geodesics under Riemannian submersion.
\end{abstract}

\maketitle

\section{Introduction}
Let $\mathbb{M}_d$, $\mathfrak{M}_d$, $\Symm_d$, $\mathbb{P}_d$, and $O_d$ denote the sets of all, invertible, symmetric, positive definite, and orthogonal $d \times d$ matrices respectively.
\subsection{The Deep Linear Network.} 
The deep linear network (DLN) is a phenomenological model of the training dynamics of neural networks \cite{arora2018convergence, arora2018optimization, cohen2023deep, menon2023geometry}. Given a depth $N \in \mathbb{Z}^+$ and width $d \in \mathbb{Z}^+$, the DLN in its simplest form consists of the ``upstairs'' (high-dimensional) parameter space of $N$ $d \times d$ matrices $\mathbb{M}_d^N$ and a ``downstairs'' (low-dimensional) observable space $\mathbb{M}_d$, linked by the product map $\phi:\mathbb{M}_d^N \to \mathbb{M}_d$ that sends
\begin{equation}\label{eq:phi}
    \mathbf{W} = (W_N, \ldots, W_1) \mapsto X = W_NW_{N-1}\cdots W_1.
\end{equation}
Indeed, $\mathbb{M}_d^N$ has $Nd^2$ degrees of freedom while $\mathbb{M}_d$ only has $d^2$, modeling overparameterization without the complications of nonlinearities and bias terms in traditional neural networks. The DLN is used to explicitly analyze the interplay of overparameterization with gradient descent of an energy function $E : \mathbb{M}_d \to \mathbb{R}$. The gradient flow in parameter space
\begin{equation}\label{eq:dln_gradient_flow}
\dot{\mathbf{W}} = -\nabla_\mathbf{W} E(X) = -\nabla_\mathbf{W} E(\phi(\mathbf{W}))
\end{equation}
can be written out explicitly via the chain rule as
\begin{equation}\label{eq:chain_rule}
\dot{W}_p = -(W_N\dots W_{p+1})^T dE(X) (W_{p-1}\dots W_1)^T
\end{equation}
for all $1 \leq p \leq N$, where $dE(X) \in \mathbb{M}_d$ is given elementwise by $dE(X)_{ij} = \frac{\partial E}{\partial X_{ij}}$. 

\subsection{Balancedness and Riemannian Geometry} Interestingly, the DLN gradient flow (Equation \ref{eq:dln_gradient_flow}) has invariant manifolds with associated metrics that illustrate the underlying geometric structure. One such invariant manifold is defined by the \textit{balancedness} constraints
\begin{equation}
    W_{p+1}^TW_{p+1} = W_pW_p^T
\end{equation}
for all $1 \leq p \leq N-1$ \cite[Theorem 1]{arora2018optimization}. When all coordinates are also rank $d$, the manifold is called the balanced manifold and is denoted as $\mathcal{M}$. Note that when $\mathbf{W} \in \mathcal{M}$, $X = \phi(\mathbf{W}) \in \mathfrak{M}_d$. 

In $\mathcal{M}$, $W_p$ have the same singular values for all $1 \leq p \leq N$ and the right singular vectors of $W_{p+1}$ align with the left singular vectors of $W_{p}$ for $1 \leq p \leq N-1$. Thus, given the SVD of a matrix $X = Q_N\Sigma Q_0^T$, we can parameterize any $\mathbf{W} \in \mathcal{M}$ such that $\phi(\mathbf{W}) = X$ using just 
\begin{equation}
    (\Lambda, Q_N, Q_{N-1}, \ldots, Q_0) \in \mathbb{R}_+^d\times O_d^{N + 1},
\end{equation}
where $\Lambda = \Sigma^\frac{1}{N}$ and $W_p = Q_{p}\Lambda Q_{p-1}^T$ are obtained by the SVD for all $1 \leq p \leq N$. Therefore, there exists a bijective map 
\begin{align}\label{eq:xi_parameterization}
\begin{split}
    \xi: \mathbb{R}^d_{+} \times O_d^{N+1} &\to \mathbb{M}_d^N, \\
    (\Lambda, Q_N, Q_{N-1}, \ldots, Q_{0}) &\mapsto (W_N, W_{N-1}, \ldots, W_1).
\end{split}
\end{align}
Because $\mathcal{M}$ is embedded within the flat $(\mathbb{M}_d^N,F)$ where $F$ is the Euclidean metric, we can induce the natural metric $\iota$ on $\mathcal{M}$ from $F$. Downstairs, we can define a metric $g^N$ on $\mathfrak{M}_d$ as follows: for any $Z \in T_X\mathfrak{M}_d$,
\begin{equation}
    \|Z\|_{g^N}^2 = \Tr\left(Z^T\mathcal{A}_{N,X}^{-1}(Z)\right),
\end{equation}
where $\mathcal{A}_{N,X}: T_X\mathfrak{M}_d^* \to T_X\mathfrak{M}_d$ is the linear operator that maps from 1-forms to the tangent space with explicit form
\begin{equation}
    Z \mapsto \sum_{p=1}^N\left(XX^T\right)^\frac{N-p}{N}Z\left(X^TX\right)^\frac{p-1}{N}.
\end{equation}
$g^N$ is crucial in the DLN because it can be used to express the training dynamics of the end-to-end matrix $X$ as a Riemannian gradient flow \cite[Section 4]{bah2022learning}. Furthermore, these metrics in tandem provide a geometric perspective on the DLN model from \cite[Theorem 8]{menon2025entropy} and \cite[Theorem 36]{tianmin2024}, which we elaborate on in Section \ref{sec:riemannian_geometry_background}:
\begin{theorem}[$\phi$ is a Riemannian Submersion]\label{thm:phi_is_submersion}
    The end-to-end product map $\phi: \mathcal{M} \to \mathfrak{M}_d$ is a Riemannian submersion from $(\mathcal{M}, \iota)$ to $(\mathfrak{M}_d, g^N)$. 
\end{theorem}

Explicit formulas for geodesics in the Bures-Wasserstein geometry can be derived by identifying a Riemannian submersion structure from invertible matrices to the positive definite matrices \cite[Section 4]{bhatia2019bures}. We present an alternate derivation of the geodesic formula that generalizes to the DLN, yielding the general ODEs for geodesics in $(\mathfrak{M}_d, g^N)$. Furthermore, Theorem \ref{thm:phi_is_submersion} implies that the DLN also has a Riemannian submersion structure with $(\mathcal{M}, \iota)$, $(\mathfrak{M}_d, g^N)$, and $\phi$, which we leverage to solve the general equations in a specific case.
\section{Main Results}
Our main results are a statement of ODEs for general geodesics in $(\mathfrak{M}_d, g^N)$ and explicit solutions when the singular vectors of the endpoints are related by a single matrix $Q \in O_d$. The solutions are obtained by finding a sufficient condition for horizontal straight lines in $\mathcal{M}$ in a manner similar to the classical theory of confocal quadrics \cite{hilbert1952geometry}.

Let $f_\alpha^{[1]}(\Lambda)$ where $\Lambda = \diag(\lambda_1, \lambda_2, \ldots, \lambda_d)$ denote the matrix
\begin{equation}\label{eq:derivative_kernel_alpha}
    \left(f_{\alpha}^{[1]}(\Lambda)\right)_{ij} = \begin{cases}
        \frac{\lambda_i^{\alpha} - \lambda_j^{\alpha}}{\lambda_i - \lambda_j} &  i \neq j, \\
        \alpha \lambda_i^{\alpha - 1} &  i = j.
    \end{cases}
\end{equation}
We discuss the origin of this matrix in Section \ref{sec:dln_general_case_proof}. Let $\circ$ denote the Hadamard product between two matrices.
\begin{theorem}[DLN Geodesic Equations]\label{thm:dln_general_case}
    $X(t) = U(t)\Sigma(t)V(t)^T$ is a geodesic between $A$ and $B$ on $(\mathfrak{M}_d, g^N)$ if and only if it satisfies the Dirichlet boundary value problem
    \begin{equation}
        \begin{cases}
            \dot{X} = \mathcal{A}_{N,X}(P), \\
            \dot{P} = -U \left(\sum_{p=1}^N \left[f_{\frac{N-p}{N}}^{[1]}\left(\Sigma^2\right) \circ M_1\right] \Sigma + \Sigma\left[f_\frac{p-1}{N}^{[1]}\left(\Sigma^2\right) \circ M_2\right]\right) V^T,        \\
            X(0) = A, \\
            X(1) = B,
        \end{cases}
    \end{equation}
    where $M_1$ and $M_2$ are the symmetric matrices
    \begin{align}\label{eq:m1_and_m2}
    M_1 &= (U^TPV)\Sigma^\frac{2(p-1)}{N}(U^TPV)^T \\
    M_2 &= (U^TPV)^T \Sigma^{\frac{2(N-p)}{N}}(U^TPV).
    \end{align}
\end{theorem}
We discover these equations via direct computation in Hamiltonian coordinates (Section \ref{sec:dln_general_case_proof}). Using Riemannian submersion, we are then able to derive explicit solutions to these equations for particular endpoints. As an important intermediate result, we find a sufficient condition for straight lines in $\mathcal{M}$.

Let $\mathbf{A} = \xi(\Lambda, Q_{N}, \ldots, Q_{0})$ and $\mathbf{B} = \xi(\tilde{\Lambda}, \tilde{Q}_N, \ldots, \tilde{Q}_0)$ be two endpoints in $\mathcal{M}$. Then, define \begin{equation}
    R_{p} = \tilde{Q}_{p}^TQ_{p} \in O_d.
\end{equation}
\begin{lemma}[Straight Lines on $\mathcal{M}$]\label{lem:characterization_of_straight_lines}

\begin{enumerate}
    \item The line \begin{equation}
        \mathbf{W}(t) = (1-t) \mathbf{A} + t \mathbf{B}
    \end{equation}
    is balanced for all $t \in [0,1]$ if and only if
    \begin{equation}
        (R_{p+1} - R_{p-1})\left(\tilde{\Lambda}^{-1} R_p \Lambda\right)^T, 
    \end{equation}
    is skew symmetric for all $1 \leq p \leq N-1$.
    \item  If $R_p = R_{p-1}$ for all $1 \leq p \leq N$, $\mathbf{W}(t) \in \mathcal{M}$ for all $t \in [0, 1]$.
\end{enumerate}
    
\end{lemma}

We are now ready to state the second main result. Let $A = U\Sigma V^T$ and $B = \tilde{U}\tilde{\Sigma}\tilde{V}^T$ be the endpoints in $\mathfrak{M}_d$ of the geodesic.
\begin{theorem}[Explicit Formulas for DLN Geodesics]\label{thm:dln_special_case_geodesics}
    If there exists $Q \in O_d$ such that $\tilde{U} = UQ$ and $\tilde{V} = VQ$,
    \begin{enumerate}
        \item there exists $\mathbf{A} \in \phi^{-1}(A)$ and $\mathbf{B} \in \phi^{-1}(B)$ such that $x(t) = (1-t)\mathbf{A} + t\mathbf{B}$ is a geodesic in $(\mathcal{M}, \iota)$.
        \item the geodesic in $(\mathfrak{M}_d, g^N)$ between $A$ and $B$ is given by 
    \begin{equation}
    \begin{split}
        X(t) = &\left((1-t)U\Sigma^\frac{1}{N} +  t\tilde{U}\tilde{\Sigma}^\frac{1}{N}Q^T\right)\\
        &\quad\left((1-t)\Sigma^\frac{1}{N} + tQ\tilde{\Sigma}^\frac{1}{N}Q^T\right)^{N-2}\\
        &\quad\left((1-t)\Sigma^\frac{1}{N}V^T+t Q\tilde{\Sigma}^\frac{1}{N}\tilde{V}^T\right).
    \end{split}
    \end{equation}
    \end{enumerate}
\end{theorem}
\begin{remark}
    Examples of when $A$ and $B$ satisfy this right-left rotation condition are if they are both symmetric or if they have the same singular vectors ($Q = \text{Id}$).
\end{remark}

\begin{remark}
When $A$ and $B$ are diagonal with positive entries (fractional powers are well defined), the geodesic simplifies to
\begin{align}
\begin{split}
    X(t) &= \left((1-t)A^{1/N}+ tB^{1/N}\right)^N \\
    &= \sum_{p=0}^N\binom{N}{p} (1-t)^pt^{N-p}A^{p/N}B^{(N-p)/N}.
\end{split}
\end{align}
This simplification can be shown by applying Theorem \ref{thm:dln_special_case_geodesics} and expanding the geodesic formulas using the commutativity of diagonal matrices.  The limiting behavior as $N \to \infty$ (recognizing $X(t)$ as the expectation of $A^SB^{1-S}$ where $S$ is the average of $N$ Bernoulli$(1-t)$ random variables) is 
\begin{equation}\label{eq:infinite_depth_geodesic}
    X(t) = A^{1-t}B^t.
\end{equation}
\end{remark}

The two important tools we will use to prove Theorem \ref{thm:dln_special_case_geodesics} in Section \ref{sec:dln_special_case_geodesics_proof} are (1) the sufficient condition for straight lines on $\mathcal{M}$ (Lemma \ref{lem:characterization_of_straight_lines}) and (2) an orthonormal basis for the horizontal space in $\mathcal{M}$ (Lemma \ref{lem:horizontal_system_of_equations} and \cite{menon2025entropy, tianmin2024}). Straight lines are automatically geodesics on $(\mathcal{M}, \iota)$ because $\iota$ is induced from the Euclidean metric on $\mathbb{M}_d^N$. We will then send a subset of these straight lines to new geodesics in $(\mathfrak{M}_d, g^N)$ to uncover the explicit formulas for geodesics as stated.

\section{Proofs}
\subsection{Geodesics Review}
For more details on geodesics, we refer the reader to \cite{do1992riemannian}. Recall the following standard definitions and results.
\begin{definition}[Geodesic]
A geodesic between $A, B \in \mathcal{X}$ on the Riemannian manifold $(\mathcal{X}, g)$ is an extrema of the action 
\begin{equation}
    \mathcal{A}[x] = \int_0^1 L(x, \dot x)   dt, \quad L(x, \dot{x}) = g_{x(t)} (\dot{x}(t), \dot{x}(t))
\end{equation}
over the set of paths
\begin{equation}
    \mathcal{S} = \{x: [0,1] \to \mathcal{X} | x(0) = A, x(1) = B, x\in \mathcal{C}^1([0,1])\}.
\end{equation}
$L$ is referred to as the \textit{Lagrangian}.
\end{definition}

\begin{theorem}[Euler-Lagrange and Hamiltonian Equations]
A geodesic $x$ satisfies the Euler-Lagrange equations:
\begin{equation}
    \frac{d}{dt}\frac{\partial L(x, \dot{x})}{\partial \dot{x}} = - \frac{\partial L(x,\dot{x})}{\partial x}.
\end{equation}
and the dual Hamiltonian equations
\begin{align}
    \dot{x} &= \partial_p H(x,p), \\
    \dot{p} &= -\partial_x H(x,p).
\end{align}
where $H$ is the Legendre transform of $L(x, \dot{x})$
\begin{equation}
    H(x, p) = \sup_{\dot{x}} (\langle p, \dot{x} \rangle - L(x, \dot{x})).
\end{equation}
with change of variables $p = \frac{\partial L}{\partial \dot{x}}$.
\end{theorem}
For a Riemannian manifold $(\mathcal{X}, g)$, the Hamiltonian $H$ associated with Lagrangian $L(x,\dot{x}) = \frac{1}{2}g(\dot{x}, \dot{x})$ is always $H(x, p) = \frac{1}{2} g^{-1}(p, p)$, where $g^{-1}$ is the inverse metric tensor acting on cotangent vectors.
\subsection{The Bures-Wasserstein Geometry}\label{sec:bw_geometry}
The Bures-Wasserstein geometry \cite{bhatia2019bures} is a special case of the DLN model when $N = 2$ and $W_1 = W_2^T$ that will serve as a motivating example for the derivation of Theorem \ref{thm:dln_general_case}. The parameter space is $\mathfrak{M}_d$ whereas the observable space is $\mathbb{P}_d$. The map $\pi: \mathfrak{M}_d \to \mathbb{P}_d$ is a Riemannian submersion that maps
\begin{equation*}
    W \in \mathfrak{M}_d \mapsto X = WW^T \in \mathbb{P}_d.
\end{equation*}
There are $d^2$ degrees of freedom in parameter space compared to only $\binom{d+1}{2}$ dimensions downstairs as $X$ is symmetric. The downstairs metric, now denoted $g^{\text{BW}}$, at a point $X \in \mathbb{P}_d$ reduces to a very simple form
\begin{equation}
g^\text{BW}(Z_1, Z_2) =\frac{1}{2}\Tr(Z_1^T \mathcal{L}^{-1}_X(Z_2)),
\end{equation}
where $\mathcal{L}^{-1}_X(Z)$ is the Lyapunov operator that maps $Z$ to the matrix that solves
\begin{equation}
    \mathcal{L}^{-1}_X(Z) X + X\mathcal{L}^{-1}_X(Z) = Z.
\end{equation}
The inverse metric tensor is $\mathcal{L}_X(Y) = XY + YX$. An important quantity in the Bures-Wasserstein geometry is the geometric mean of two matrices.
\begin{definition}[Geometric Mean]\label{def:geometric_mean}
Given $A, B \in \mathbb{P}_d$, the geometric mean of $A$ and $B$ (denoted by $A \# B$) is
\begin{equation}
    A\#B = A^{1/2}(A^{-1/2}BA^{-1/2})^{1/2}A^{1/2}.
\end{equation}
\end{definition}
We will use two key facts about the geometric mean, both of which can be proven directly from Definition \ref{def:geometric_mean}.
\begin{lemma}\label{lem:geometric_mean_riccati}
\begin{enumerate}
    \item $A\#B$ is the unique positive definite solution $M$ to the Riccati Equation
    \begin{equation}
        B = MA^{-1}M.
    \end{equation}
    \item $A(A^{-1}\#B) = (AB)^{1/2}$ and $(A^{-1}\#B)A = (BA)^{1/2}$.
\end{enumerate}
\end{lemma}
Our main computation in the Bures-Wasserstein geometry is a derivation of the explicit formulas for geodesics by solving the Hamiltonian equations of motion directly. 
\begin{theorem}[Bures-Wasserstein Geodesics]\label{thm:bw_general_case}
    The geodesic path $X: [0,1] \to \mathbb{P}_d$ between $A, B \in \mathbb{P}_d$ in $(\mathbb{P}_d, g^\text{BW})$ satisfies the boundary value problem
    \begin{equation}\label{eq:bures_wasserstein_equations}
        \begin{cases}
            \dot{X} = \mathcal{L}_X(P), \\
            \dot{P} = -P^2,\\
            X(0) = A,\\
            X(1) = B,
        \end{cases}
    \end{equation}
    and has explicit form given by
\begin{equation}
    X(t) = (1-t)^2 A + t^2 B + t(1-t)\left((AB)^{1/2} + (BA)^{1/2}\right).
\end{equation}
\end{theorem}

\begin{proof}
To begin, recall that $T_X\mathbb{P}_d \cong \text{Symm}_d$ and denote $Y = \dot{X}$. We consider the Hamiltonian system, which from the inverse metric tensor $\mathcal{L}$ is
\begin{align}
    H(X, P) &= \frac{1}{2}\Tr(P\mathcal{L}_X(P)) = \Tr(P^2X).
\end{align}
The system of differential equations is now clear from the form of $H(X,P)$. We can solve them to find exact formulas for $X$ and $P$. In particular, $\dot{P} = -P^2$ implies
\begin{align}
    -P^{-1} \dot{P} P^{-1} &= I \nonumber \\
    \frac{d}{dt}(P^{-1}) &= I \nonumber \\
    P^{-1}(t) &= tI + P^{-1}_0.
\end{align}
Thus, $P(t) = (tI + P^{-1}_0)^{-1}$. Factoring out $P^{-1}_0$ from the interior term reveals that $P(t)$ is a time derivative:
\begin{align}
    P(t) &= (P^{-1}_0(tP_0 + I))^{-1}  \nonumber \\
    &= (tP_0 + I)^{-1}P_0 \\
    &= \frac{d}{dt} \log (tP_0 + I).
\end{align}
The first differential equation can be expanded as
\begin{align}
    \dot{X} = X(t)P(t) + P(t)X(t).
\end{align}
Motivated by the similarity of this equation to the Lax equation, we make the ansatz for the solution form of $X$ as
\begin{equation}
    X(t) = e^{C(t)} X_0 e^{D(t)},
\end{equation}
which after plugging in and matching coefficients gives that $\dot{C}(t) = \dot{D}(t) = P(t)$ and implying that 
\begin{align}
    C(t) = D(t) &= \int_0^t P(s) ds\\
    &= \int_{0}^t \frac{d}{ds}\log(sP_0 + I) ds \\
    &= \log (tP_0+I).
\end{align}
Substituting $C$ and $D$ back into the ansatz gives that 
\begin{align}
    X(t) &= \exp\{\log (tP_0+I)\} X_0 \exp\{\log (tP_0+I)\}\nonumber \\
    &= (tP_0 + I) X_0 (tP_0 + I) \nonumber \\ 
    &= (tP_0 + I) A (tP_0 + I).
\end{align}
Applying the final unused condition that $X_1 = B$ helps find $P_0$:
\begin{equation}
    B = (P_0 + I) A (P_0 + I),
\end{equation}
which from Lemma \ref{lem:geometric_mean_riccati} implies that
\begin{equation}
    P_0 = (A^{-1} \# B) - I.
\end{equation}
Finally, substituting $P_0$ into our solution formula for $X(t)$, we conclude that
\begin{align}
    X(t) &= (t((A^{-1} \# B) - I) + I) A (t((A^{-1} \# B) - I) + I) \\
    &= (t(A^{-1} \# B) + (1-t)I) A (t(A^{-1} \# B) + (1-t)I)\\
    &= (1-t)^2 A + t^2 B + t(1-t)((A^{-1}\#B)A + A(A^{-1}\#B))\\
    &= (1-t)^2 A + t^2 B + t(1-t)\left((AB)^{1/2} + (BA)^{1/2}\right).
\end{align}
The uniqueness of this solution also follows from Picard's theorem, as the right hand sides of Equation \ref{eq:bures_wasserstein_equations} are smooth.
\end{proof}

\subsection{General Equations from Hamiltonian}\label{sec:dln_general_case_proof}
We will now generalize the preceding calculation to the DLN. First, recall functional calculus on symmetric matrices \cite[Chapter 2]{bhatia2009positive}. If $S = Q\Lambda Q^T$ is a symmetric matrix and $f : \mathbb{R} \to \mathbb{R}$ is a 1D functional, $f$ can be extended to the symmetric matrices as (abusing notation)
    \begin{equation}
    f(S) = Qf(\Lambda)Q^T = Q\text{diag}(f(\lambda_1), f(\lambda_2), \dots, f(\lambda_d))Q^T.
    \end{equation}
    Then, the differential of $f(S)$ at $S$ evaluated in the direction $Z$ is exactly
    \begin{equation}
        Df[S](Z) = Q\left(f^{[1]}(\Lambda) \circ \left[Q^T Z Q\right]\right) Q^T
    \end{equation}
    where $f^{[1]}$ is the special matrix given elementwise by the conditional
    \begin{equation}\label{eq:derivative_kernel}
        f^{[1]}(\Lambda)_{ij} = \begin{cases}
            \frac{f(\lambda_i) - f(\lambda_j)}{\lambda_i - \lambda_j} & i \neq j \\
            \lim_{\lambda_i \to \lambda_j}\frac{f(\lambda_i) - f(\lambda_j)}{\lambda_i - \lambda_j} & i = j.
        \end{cases}
    \end{equation}
Equation \ref{eq:derivative_kernel_alpha} for $f^{[1]}_\alpha$ is the special case of Equation \ref{eq:derivative_kernel} when $f(x) = x^\alpha$. We also have the following derivative computation.
\begin{lemma}[Derivatives of $(XX^T)^\alpha$]\label{lem:deriv_of_alpha}
Let $X = U\Sigma V^T$ be the singular value decomposition of $X$. Then,
    \begin{align}
D_X\left[\left(XX^T\right)^{\alpha}\right](Z) &= U\left[f_{\alpha}^{[1]}\left(\Sigma^2\right) \circ \left(\Sigma \tilde{Z}^T + \tilde{Z}\Sigma\right)\right] U^T.\\        D_X\left[\left(X^TX\right)^{\alpha}\right](Z) &= V\left[f_{\alpha}^{[1]}\left(\Sigma^2\right) \circ \left(\Sigma \tilde{Z} + \tilde{Z}^T \Sigma\right)\right] V^T.
    \end{align}
    where $\tilde{Z} = U^TZV$.
\end{lemma}
\begin{proof}
    We write out the proof for the first result - the second is similar. Notice that $XX^T$ is a symmetric matrix with diagonalization $U\Sigma^2U^T$. Thus, we apply functional calculus along with the chain rule to find that 
    \begin{equation}
        D_X\left[\left(XX^T\right)^\alpha\right](Z) = U\left[f_\alpha^{[1]}\left(\Sigma^2\right) \circ \left( U^T\left(D_X\left[XX^T\right](Z) \right)U\right)\right] U^T.
    \end{equation}
    Clearly, we know that 
    \begin{equation}
        D_X\left[XX^T\right](Z) = XZ^T + ZX^T.
    \end{equation}
    The final simplification comes from expanding using $Z = U \tilde{Z} V^T$.
\end{proof}

We also prepare another intermediate result that will be useful for algebraic manipulation.
\begin{lemma}\label{lem:some_trace_identities}
    Let $K$ and $M$ be symmetric matrices, $\Sigma$ be a diagonal matrix, and $A$ be an arbitrary matrix. Then, we have the following identities.
    \begin{align}
        \Tr([K \circ (\Sigma A^T)]M) &= \Tr([K \circ (A\Sigma)]M) = \Tr(\Sigma[K \circ M] A)
    \end{align}
\end{lemma}
\begin{proof}
    First, we prove that the first expression is equal to the last. We will approach by expanding into coordinates.
    \begin{align}
        \Tr([K \circ (\Sigma A^T)] M) &= \sum_{i} \sum_j K_{ij}\sigma_i(A^T)_{ij}M_{ji}\nonumber\\
        &= \sum_{i} \sum_{j} \sigma_iK_{ij}M_{ji} A_{ji} \nonumber\\
        &= \sum_i \sum_j \sigma_i K_{ij} M_{ij} A_{ji} \nonumber\\
        &= \Tr(\Sigma [K\circ M] A).
    \end{align}
    Similarly, we can prove the equality of the second expression and the third expression using the symmetry of $K$.
    \begin{align}
        \Tr([K \circ (A\Sigma)] M) &= \sum_{i} \sum_j K_{ij}A_{ij}\sigma_jM_{ji}\nonumber\\
        &= \sum_i \sum_j \sigma_jK_{ij}M_{ji} A_{ij} \nonumber
        \\
        &= \sum_i \sum_j \sigma_j K_{ji} M_{ji} A_{ij} \nonumber \\
        &= \Tr(\Sigma [K\circ M] A).
    \end{align}
\end{proof}

The Hamiltonian on the manifold $(\mathfrak{M}_d, g^N)$ is again given by the inverse metric tensor:
\begin{equation}
    H(X,P) = \frac{1}{2}\Tr\left(P^T \mathcal{A}_{N,X}(P)\right).
\end{equation}
Using these coordinates, we will now obtain the general equations on the manifold $(\mathfrak{M}_d, g^N)$ through direct computation, proving Theorem \ref{thm:dln_general_case}.
\begin{proof}[Proof of Theorem \ref{thm:dln_general_case}]
The first equation follows directly from the Hamiltonian. The second equation requires computation of $D_X[H(X,P)](Z)$. We begin by applying Lemma \ref{lem:deriv_of_alpha}.
\begin{align}
    &D_X[H(X,P)](Z) = \frac{1}{2}\sum_{p=1}^N \Tr\left(T_U + T_V\right) = \frac{1}{2}\sum_{p=1}^N \Tr(T_U) + \Tr(T_V),
    \\&T_U = P^TU\left[f_{\frac{N-p}{N}}^{[1]} \left(\Sigma^2\right) \circ \left(\Sigma \tilde{Z}^T + \tilde{Z}\Sigma\right)\right]U^T P V\Sigma^{\frac{2(p-1)}{N}}V^T,
    \\&T_V = P^TU\Sigma^{\frac{2(N-p)}{N}}U^T P V\left[f_{\frac{p-1}{N}}^{[1]}\left(\Sigma^2\right) \circ \left(\Sigma \tilde{Z} + \tilde{Z}^T \Sigma\right)\right]V^T.
\end{align}
To clean up further, we can write $P$ in terms of the singular value coordinates by setting $\tilde{P} = U^T P V$ and using the cyclic property of trace to see that
\begin{align}
    \Tr(T_U) &= \Tr\left(\tilde{P}^T \left[f_{\frac{N-p}{N}}^{[1]}\left(\Sigma^2\right) \circ \left(\Sigma \tilde{Z}^T + \tilde{Z}\Sigma\right)\right] \tilde{P}\Sigma^{\frac{2(p-1)}{N}}\right)\label{eq:dln_geodesic_inter1}\\
    \Tr(T_V) &= \Tr\left(\tilde{P}^T \Sigma^{\frac{2(N-p)}{N}} \tilde{P} \left[ f_{\frac{p-1}{N}}^{[1]} \left(\Sigma^2\right) \circ \left(\Sigma \tilde{Z} + \tilde{Z}^T \Sigma\right)\right]\right)\label{eq:dln_geodesic_inter2}.
\end{align}

Applying the identities from Lemma \ref{lem:some_trace_identities} to Equations \ref{eq:dln_geodesic_inter1} and \ref{eq:dln_geodesic_inter2} ($f^{[1]}_{\alpha}(\Sigma^2)$ is symmetric) and substituting in $M_1$ and $M_2$ from Equation \ref{eq:m1_and_m2} shows
\begin{align}
    D_X[H(X,P)](Z) &= \frac{1}{2}\sum_{p=1}^N \Tr\left(2\Sigma \left[f_{\frac{N-p}{N}}^{[1]}\left(\Sigma^2\right) \circ M_1\right]\tilde{Z} + 2\left[f_\frac{p-1}{N}^{[1]}\left(\Sigma^2\right) \circ M_2\right]\Sigma\tilde{Z}\right) \\
    &= \sum_{p=1}^N \Tr\left(\left(\Sigma \left[f_{\frac{N-p}{N}}^{[1]}\left(\Sigma^2\right) \circ M_1\right] + \left[f_\frac{p-1}{N}^{[1]}\left(\Sigma^2\right) \circ M_2\right]\Sigma\right)\tilde{Z}\right)\\
    &= \Tr\left(V\sum_{p=1}^N \left(\Sigma \left[f_{\frac{N-p}{N}}^{[1]}\left(\Sigma^2\right) \circ M_1\right] + \left[f_\frac{p-1}{N}^{[1]}\left(\Sigma^2\right) \circ M_2\right]\Sigma\right)U^T Z\right).
\end{align}
This lets us conclude that the second canonical equation in the Hamiltonian formulation is 
\begin{equation}
    \dot{P} = -U \sum_{p=1}^N \left(\left[f_{\frac{N-p}{N}}^{[1]}\left(\Sigma^2\right) \circ M_1\right] \Sigma + \Sigma\left[f_\frac{p-1}{N}^{[1]}\left(\Sigma^2\right) \circ M_2\right]\right) V^T.
\end{equation}
\end{proof}
However, unlike the Bures-Wasserstein case (Theorem \ref{thm:bw_general_case}), we find that this final form is not particularly amenable to direct solution, despite the commutator structure and relation to the singular vector coordinates $U$ and $V$. Thus, we turn to special cases, where we will be able to find explicit formulas by leveraging the inherent geometric structure of the DLN model.

\subsection{Explicit Solutions from Submersion}\label{sec:dln_special_case_geodesics_proof}
\subsubsection{Riemannian Geometry Background}\label{sec:riemannian_geometry_background}
Recall the following important definitions and results from Riemannian geometry \cite{bhatia2019bures}.
\begin{definition}[Vertical and Horizontal Spaces]
 Let $(\mathcal{X}, g)$ and ($\mathcal{Y}, h)$ be linked by a differentiable map $\pi: \mathcal{X} \to \mathcal{Y}$. Then, there exists a direct sum decomposition of the tangent space at any point $m \in \mathcal{X}$ via the differential of $\phi$.
\begin{equation}
    T_m \mathcal{X} = \ker D\pi[m] \oplus (\ker D\pi[m])^\perp = \mathcal{V}_m \oplus \mathcal{H}_m,
\end{equation}
where $\mathcal{V}_m$ and $\mathcal{H}_m$ denote the vertical and horizontal spaces at $m$ respectively.
\end{definition}
\begin{definition}[Riemannian Submersion]
Let $(\mathcal{X}, g)$ and $(\mathcal{Y},h)$ be two Riemannian manifolds and let $T_m\mathcal{X} = \mathcal{V}_m \oplus \mathcal{H}_m$ be the direct sum decomposition of the tangent space at $m \in \mathcal{X}$. A mapping $\pi : \mathcal{X} \to \mathcal{Y}$ is a Riemannian submersion if and only if the differential $D\pi[m]: T_m\mathcal{X} \to T_{\pi(m)}\mathcal{Y}$ satisfies the following for all $m \in \mathcal{X}$:
\begin{enumerate}
    \item $D\pi[m]$ is surjective ($\pi$ is a smooth submersion).
    \item $D\pi[m]$ restricted to the horizontal space is an isometry i.e.\ $\|z\|_g = \|D\pi[m](z)\|_h$ for all $z \in \mathcal{H}_m$.
\end{enumerate}
\end{definition}

A core fact underlying \cite[Section 4]{bhatia2019bures} and Theorem \ref{thm:dln_special_case_geodesics} is that understanding the horizontal space of a submersion allows us to identify geodesics that are preserved under the submersion, as described by the following theorem \cite{oneill1966fundamental, bhatia2019bures}.
\begin{theorem}[Submersions Preserve Horizontal Geodesics]\label{thm:horizontal_geodesic_submersion}
    Let $(\mathcal{X}, g)$ and $(\mathcal{Y}, h)$ be Riemannian manifolds and $\phi: \mathcal{X} \to \mathcal{Y}$ be a Riemannian submersion. Let $x: [0,1] \to \mathcal{X}$ be a geodesic between $A, B \in \mathcal{X}$ in the geometry $(\mathcal{X}, g)$ so that $x'(0) \in \mathcal{H}_{A}$ i.e.\ is horizontal. Then,
    \begin{enumerate}
        \item $x'(t)$ is horizontal for all $t \in [0,1]$.
        \item $\phi \circ x$ is a geodesic in $(\mathcal{Y}, h)$ between $\phi(A)$ and $\phi(B)$ and is the same length as $x$.
    \end{enumerate}
\end{theorem}
Thus, if we are able to classify horizontal geodesics in $(\mathcal{M}, \iota)$, Theorem \ref{thm:horizontal_geodesic_submersion} will give explicit formulas for geodesics in $(\mathfrak{M}_d, g^N)$ through the submersion $\phi$. The simplest class of geodesics in $(\mathcal{M}, \iota)$ are straight lines. Thus, we begin with a proof of Lemma \ref{lem:characterization_of_straight_lines} that describes a set of straight lines in $\mathcal{M}$.
\begin{proof}[Proof of Lemma \ref{lem:characterization_of_straight_lines}]
    First, we wish to find the equivalent condition to $\mathbf{W}(t) = (1-t)\mathbf{A} + t\mathbf{B}$ being balanced for all $t \in (0,1)$ when $\mathbf{A}, \mathbf{B} \in \mathcal{M}$. Fix $p \in \{1, 2, \ldots, N-1\}$. We would like to find an equivalent condition to
    \begin{equation}
        W_{p+1}(t)^TW_{p+1}(t) = W_p(t)W_p(t)^T,
    \end{equation}
    where $W_p(t) = (1-t)W_p(0) + tW_p(1)$ for all $p$. Expanding the algebra out and using that $W_p(0) = A_p$ and $W_p(1) = B_p$ are both in the balanced manifold, we find that $\mathbf{W}(t)$ is balanced for all $t \in (0,1)$ if and only if
    \begin{equation}
        B_{p+1}^TA_{p+1} + A_{p+1}^TB_{p+1} = B_{p}A_{p}^T + A_{p}B_{p}^T.
    \end{equation}
    Now, we use the parameterizations of $A$ and $B$ to expand into the orthogonal matrix components.
    \begin{align}
    \begin{split}
        \tilde{Q}_{p}\tilde{\Lambda}(\tilde{Q}_{p+1}^T Q_{p+1})\Lambda Q_{p}^T &+ Q_{p}\Lambda (Q_{p+1}^T\tilde{Q}_{p+1})\tilde{\Lambda}\tilde{Q}_{p}^T \\
        &= 
        \tilde{Q}_{p}\tilde{\Lambda}(\tilde{Q}_{p-1}^T Q_{p-1})\Lambda Q_{p}^T + Q_{p}\Lambda (Q_{p-1}^T\tilde{Q}_{p-1})\tilde{\Lambda}\tilde{Q}_{p}^T.
    \end{split}
    \end{align}
    Rewriting using $R_p$ notation and left and right multiplying by $\tilde{Q}_p^T$ and $Q_p$ respectively gives 
    \begin{align}
        \tilde{\Lambda}(
R_{p+1} - R_{p-1})\Lambda  &= - R_p \Lambda (R_{p+1} - R_{p-1})^T \tilde{\Lambda} R_p.
    \end{align}
    We can rewrite in terms of $M = \tilde{\Lambda}^{-1} R_p \Lambda$ to get the final equivalent condition:
    \begin{align}
        R_{p+1} - R_{p-1} &= -M(R_{p+1} - R_{p-1})^TM^{-T}.\\
        (R_{p+1} - R_{p-1})M^T &= -M(R_{p+1} - R_{p-1})^T.
    \end{align}
    The structure demonstrates that $\mathbf{W}(t)$ is balanced if and only if the difference 
    \begin{equation}
        (R_{p+1} - R_{p-1})M^T
    \end{equation}
    is a skew symmetric matrix.

    It remains to show that if $R_p = R_{p-1}$ for all $1 \leq p \leq N$, $\mathbf{W}(t) \in \mathcal{M}$ for all $t \in [0, 1]$. In this case, the above skew symmetry condition is immediate as all differences $R_{p+1} - R_{p-1}$ are $0$, implying that $\mathbf{W}(t)$ is balanced for all $t \in [0,1]$. Thus, we just need to show that the coordinates of $\mathbf{W}(t)$ stay full rank at all intermediate $t$.

    Fix $1 \leq p \leq N$ and consider the $p$-th coordinate of $\mathbf{W}(t)$. Since $R_p$ is the same for all $p$, we drop the index for clarity.
    \begin{align}
        \mathbf{W}_p(t) &= (1-t) A_p + tB_p \\
        &= (1-t) Q_p\Lambda Q_{p-1}^T + t\tilde{Q}_p\tilde{\Lambda}\tilde{Q}_{p-1}^T \\
        &= (1-t) \tilde{Q}_p R \Lambda R^T \tilde{Q}_{p-1}^T + t\tilde{Q}_p \tilde{\Lambda} \tilde{Q}_{p-1}^T \\
        &= \tilde{Q}_p \left((1-t) R\Lambda R^T + t \tilde{\Lambda}\right)\tilde{Q}_{p-1}^T.
    \end{align}
    Because $\mathbb{P}_d$ is convex, the interior convex combination is also in $\mathbb{P}_d$ as $R\Lambda R^T \in \mathbb{P}_d$. Since $\tilde{Q}_p, \tilde{Q}_{p-1} \in O_d$ and are invertible, we conclude that $\mathbf{W}_p(t)$ is full rank for all $t \in [0,1]$ and $\mathbf{W}(t) \in \mathcal{M}$ follows.
\end{proof}

We now state a basis for the horizontal space in $\mathcal{M}$ computed in \cite[Section 6]{menon2025entropy} and \cite[Theorem 34]{tianmin2024}. Let 
\begin{equation}
    \mathbf{M} \in \mathcal{M}, \mathbf{M} = \left(Q_N\Lambda^\frac{1}{N} Q_{N-1}^T, Q_{N-1}\Lambda^\frac{1}{N} Q_{N-2}^T, \dots, Q_1\Lambda^\frac{1}{N} Q_0^T\right)
\end{equation}
so that $\mathbf{M} = \xi(\Lambda, Q_N, Q_{N-1}, \ldots, Q_1)$ and $\Lambda = \diag(\lambda_1, \lambda_2, \ldots, \lambda_d)$. The horizontal space at $\mathbf{M}$ is $\mathcal{H}_\mathbf{M} = \ker (\phi^* [\mathbf{M}])^\perp$. Let $q_{p,i}$ denote the $i$-th column of $Q_p$.  
\begin{lemma}[Basis for Horizontal Space in $\mathcal{M}$]\label{lem:horizontal_system_of_equations}
Let $Z$ be an arbitrary matrix in $\mathbb{M}_d^N$. Then, $Z \in \mathcal{H}_\mathbf{M}$ if and only if there exists a set of coefficients $\alpha \in \mathbb{M}_d$ such that
\begin{equation}\label{eq:horizontal_system_of_equations}
    Q_p\alpha Q_{p-1}^T = Z_p
\end{equation}
for all $1 \leq p \leq N$. Equivalently, $\mathcal{H}_\mathbf{M}$ is spanned by
\begin{equation}
    \left\{\left(q_{N,k} q_{N-1, l}^T, q_{N-1,k} q_{N-2, l}^T, \ldots q_{1,k} q_{0, l}^T\right) | 1 \leq k, l \leq d \right\}.
\end{equation}
\end{lemma}

We can now use $\phi$ and Theorem \ref{thm:horizontal_geodesic_submersion} to construct new explicit formulas for geodesics on $\mathfrak{M}_d$.
\begin{proof}[Proof of Theorem \ref{thm:dln_special_case_geodesics}]
    It is sufficient to first construct $\mathbf{A} \in \phi^{-1}(A) \subseteq \mathcal{M}$ and $\mathbf{B} \in \phi^{-1}(B) \subseteq \mathcal{M}$ such that \begin{equation}
    x(t) = (1-t)\mathbf{A}+ t\mathbf{B} \in \mathcal{M} \text{ for all } t \in [0,1]
    \end{equation}
    because $\iota$ is induced by the flat Frobenius metric. Then, we will show that $x'(t)$ is horizontal at $t=0$ by proving that
    \begin{equation}\label{eq:horizontal_derivative_condition}
        \frac{d}{dt} \left((1-t) \mathbf{A} + t\mathbf{B}\right) = -\mathbf{A} + \mathbf{B} \in \mathcal{H}_{\mathbf{A}}
    \end{equation}
    or that equivalently $\mathbf{A} \in \mathcal{H}_{\mathbf{A}}$ and $\mathbf{B} \in \mathcal{H}_{\mathbf{A}}$ because $\phi^*$ is linear. Finally, Theorem \ref{thm:horizontal_geodesic_submersion} can be applied, implying that $\phi \circ x$ is a geodesic in the downstairs space between $\phi(\mathbf{A})$ and $\phi(\mathbf{B})$ and giving the explicit formula as desired.
    
    We begin by identifying the correct parameterizations for the upstairs endpoints. Consider the following parameterization for $\mathbf{A}$:
    \begin{equation}\label{eq:underline_A_parameterization}
        \mathbf{A} = \xi \left(\Sigma^\frac{1}{N}, U, I, \ldots, I, V\right)
    \end{equation}
    where $\xi$ is the map given in Equation \ref{eq:xi_parameterization} since $\mathbf{A} \in \mathcal{M}$. Clearly, $\phi(\mathbf{A}) = A$. 
    Furthermore, by Lemma \ref{lem:horizontal_system_of_equations}, $\mathbf{Z} \in \mathbb{M}_d^N$ is in $\mathcal{H}_\mathbf{A}$ if and only if there exists $\alpha \in \mathbb{M}_d$ such that
    \begin{align}
        U\alpha &= Z_N,\\
        \alpha &= Z_p, 2 \leq p \leq N-2, \\
        \alpha V^T &= Z_1.
    \end{align}
    It is immediate that $\mathbf{A}$ is in its own horizontal space $\mathcal{H}_\mathbf{A}$ using the coefficient matrix $\alpha = \Sigma^\frac{1}{N}$.
    
    Now, also consider \begin{equation}\label{eq:underline_B_parameterization}
        \mathbf{B} = \xi\left(\tilde{\Sigma}^\frac{1}{N}, \tilde{U}, Q, \ldots, Q, \tilde{V}\right),
    \end{equation}
    where $Q = U^T\tilde{U} = V^T\tilde{V}$ by assumption. Clearly, $\phi(\mathbf{B}) = B$. Furthermore, $\mathbf{A}$ and $\mathbf{B}$ satisfy the sufficient condition of Lemma \ref{lem:characterization_of_straight_lines} since $Q^T = \tilde{U}^TU = \tilde{V}^TV$ and $R_p$ is constant for all $p$, so we conclude that $x(t) \in \mathcal{M}$ for all $t \in [0, 1]$. Thus, the line $(1-t)\mathbf{A} + t\mathbf{B}$ is indeed the geodesic between $\mathbf{A}$ and $\mathbf{B}$ in $(\mathcal{M}, \iota)$ because $\iota$ is flat.

    Now we will show $x'(0)$ is horizontal for the chosen $\mathbf{A}$ and $\mathbf{B}$. We have already shown that $\mathbf{A}$ is in $\mathcal{H}_\mathbf{A}$. By Equation \ref{eq:horizontal_derivative_condition}, it is enough to show that $\mathbf{B} \in \mathcal{H}_\mathbf{A}$ as well. Again, Lemma \ref{lem:horizontal_system_of_equations} provides a sufficient system of equations to solve with $\mathbf{Z} = \mathbf{B}$:
    \begin{align}
        U\alpha &= \tilde{U}\tilde{\Sigma}^\frac{1}{N}Q^T,\\
        \alpha &= Q \tilde{\Sigma}^\frac{1}{N} Q^T, \\
        \alpha V^T &= Q \tilde{\Sigma}^\frac{1}{N} \tilde{V}^T.
    \end{align}
    The solution given by the middle equations of $\alpha = Q\tilde{\Sigma}^\frac{1}{N}Q^T$ can be checked as consistent with the others by definition of $Q$. Thus, $\mathbf{B} \in \mathcal{H}_{\mathbf{A}}$.
    
    Since we have shown that both $\mathbf{B}$ and $\mathbf{A}$ are in $\mathcal{H}_{\mathbf{A}}$ and that the geodesic linking them is a straight line, it follows from Equation \ref{eq:horizontal_derivative_condition} that $x'(0)$ is horizontal. 
    Furthermore, note that because all coordinates in $x(t) \in \mathcal{M}$, $\phi \circ x \in \mathfrak{M}_d$ for all $t \in [0, 1]$ as $\phi \circ x$ becomes a product of invertible matrices. Applying Theorem \ref{thm:horizontal_geodesic_submersion}, we conclude that $\phi \circ x$ is a geodesic in $(\mathfrak{M}_d, g^N)$ of the same length as $x$ and the explicit formula comes directly from applying $\phi$ to $x$ between our parameterized $\mathbf{A}$ and $\mathbf{B}$.
\end{proof}
\section{Open Problems} 
The main follow up work is (1) solving  the general case geodesic equations for arbitrary endpoints and (2) generalizing geodesic results to low rank (rank $< d$) matrices, not just $\mathcal{M}$, as the Riemannian metric and submersion structure exists in all ranks \cite{bah2022learning, menon2025entropy, tianmin2024}. Once geodesics and distances are understood in the geometries underlying the DLN, one can consider reconnecting them to the vector field given by the gradient flow to make conclusions about training dynamics, similar to \cite[Section 3.3]{karmarkar1990riemannian}.
\section*{Acknowledgements}
This work is based on a chapter of the author’s undergraduate thesis supervised by Govind Menon (Brown University), to whom the author is grateful for the initial problem statement and continual feedback on the project and paper. The author is also thankful for advice from Austin Stromme (ENSAE/CREST) and Kathryn Lindsey (Boston College) at various stages of this work. This work was partially supported by NSF grant 2407055.
\bibliographystyle{amsplain}
\bibliography{refs}
\end{document}